\documentclass[11pt]{amsart}

\usepackage{hyperref}
\usepackage{amsthm} 
\usepackage{amssymb}
\usepackage{graphicx}				 
\usepackage{stmaryrd}                  
\usepackage{footmisc}                  
\usepackage{paralist}
\usepackage{wrapfig}
\usepackage[draft]{pdfcomment}
\usepackage{bbm}
\usepackage[T1]{fontenc}
\usepackage[utf8]{inputenc}
\usepackage[ngerman,british]{babel}
\usepackage[arrow, matrix, curve]{xy}
\usepackage{color}
\usepackage{mathtools}
\usepackage[pagewise]{lineno}

\usepackage{enumitem}

\newtheorem{thm}{Theorem}[section]

\newtheorem{lem}[thm]{Lemma}

\newtheorem{cor}[thm]{Corollary}

\theoremstyle{definition}

\theoremstyle{remark}

\newtheorem{qst}[thm]{Question}

\newcommand{\R}{\mathbb{R}}
\newcommand{\N}{\mathbb{N}}
\newcommand{\Z}{\mathbb{Z}}
\newcommand{\Q}{\mathbb{Q}}

\newcommand{\Orb}{\mathcal{O}}

\newcommand{\To}{\rightarrow}

\newcommand{\MTo}{\mapsto}

\title[Closed geodesics on compact orbifolds and on noncompact manifolds]{Closed geodesics on compact orbifolds \\ and on noncompact manifolds}
\author[Christian Lange]{Christian Lange}
\address[Christian Lange]{Ludwig-Maximilians-Universit\"at M\"unchen, Mathematisches Institut\newline\indent Theresienstr. 39, 80333 München, Germany}
\email{lange@math.lmu.de}
\author[Christoph Zwickler]{Christoph Zwickler}
\address[Christoph Zwickler]{Universit\"at zu K\"oln, Mathematisches Institut, \newline\indent Weyertal 86-90, 50931 K\"oln, Germany}

\begin{document}

\begin{abstract} We prove the existence of a nonconstant closed geodesic on every odd-dimensional compact Riemannian orbifold and on every connected, noncontractible Riemannian manifold that admits a cocompact, isometric group action. We deduce new restrictions on a compact Riemannian orbifold without such a geodesic; in particular, its orbifold fundamental group is a rational Poincaré duality group.
\end{abstract}

\maketitle

\section{Introduction}

The study of closed geodesics has been a central problem in Riemannian geometry since the pioneering works of Hadamard \cite{Ha98} and Poincaré \cite{Poincare}. By the classical Lyusternik--Fet theorem every closed Riemannian manifold has a closed geodesic \cite{LuFet}; closed geodesics being assumed nonconstant throughout. A possible generalization of this problem is to look for closed  geodesics on Riemannian orbifolds \cite{Borzellino,Borzellino2}. Most notably, Guruprasad and Haefliger \cite{GH06} showed the existence of a closed geodesic on every nondevelopable compact Riemannian orbifold, i.e. one which cannot be written as a global quotient $M/G$ of a Riemannian manifold $M$ by a proper, cocompact, isometric action of a discrete group $G$; see also \cite{Asselle}. 

In the case of a developable compact Riemannian orbifold $M/G$ as above the existence of a closed geodesic on $M/G$ is equivalent to the existence of a nonconstant geodesic $c$ on $M$ and an isometry $g\in G$ with $g(c(t))=c(t+1)$ for all $t\in \R$; cf. Subsection \ref{sec:orbifold_geodesics} and see \cite{GT78,Gr74,Gr85, Ma13, Ma15} for existence results on such \emph{isometry-invariant} geodesics. For instance, if $M$ is odd-dimensional and has a finite-dimensional rational cohomology ring, then every ``proper'' isometry of $M$, i.e. one whose displacement function is proper, has an invariant geodesic by \cite[Theorem~3.4]{Gr85}. In particular, under these assumptions, if $G$ contains a ``proper'' isometry, then the Riemannian orbifold $M/G$ has a closed geodesic. Moreover, by a result of Dragomir every compact Riemannian orbifold of dimension $3$, $5$ or $7$ has a closed geodesic \cite{Dragomir2,Dragomir}. We advance this line of research by generalizing Dragomir's result to all odd dimensions.

\begin{thm}\label{thm:A} Every odd-dimensional compact Riemannian orbifold has a closed geodesic. 
\end{thm}

In comparison with \cite[Theorem~3.4]{Gr85} 
note that, for instance, a translation of $\R$,  a reflection of $\R^2$ and a rotation of the surface of revolution obtained by rotating the graph of $1/(1-x^2)$, $x\in (-1,1)$, around the $y$-axis are not proper in this sense. Here, the last example differs from the first two in that the isometry has an isolated fixed point; cf. the discussion below and Question \ref{qst:existence}.

We prove Theorem \ref{thm:A} by induction on the dimension with the base case being the existence of a closed geodesic on a compact Riemannian manifold of dimension at least one; see Section~\ref{sec:closed_geodesics}. The parity assumption in the theorem is used to rule out the case of a zero-dimensional Riemannian manifold, in which all curves are constant. The proof of Theorem \ref{thm:A}  moreover reduces the existence problem of closed geodesics on compact Riemannian orbifolds to the case of even-dimensional compact Riemannian orbifolds with isolated singularities. This reduction was, however, already observed by Dragomir \cite[p.~137]{Dragomir2}. 

A natural approach to prove the existence of a closed geodesic on a developable Riemannian orbifold $M/G$ is to look for closed geodesics on the Riemannian manifold $M$. Several authors have established conditions that guarantee the existence of a closed geodesic on a noncompact Riemannian manifold \cite{Benci,Secchi,Thorb1,Thorb}. The next result applies beyond the closed geodesic problem on Riemannian orbifolds, since we assume neither properness of the action nor discreteness of 
$G$.

\begin{thm}\label{thm:B} Suppose a group $G$ acts isometrically and cocompactly on a connected Riemannian manifold $M$. If $M$ is not contractible, then $M$ has a closed geodesic.
\end{thm}

As a consequence of Theorem \ref{thm:A}, Theorem \ref{thm:B}, \cite{GH06} and Lemma \ref{lem:pd_group} below, we obtain the following corollary; see Subsection \ref{sub:proof_cor}.

\begin{cor}\label{cor:poincare_duality} If a compact, connected Riemannian $n$-orbifold does not have a closed geodesic, then it is developable, its universal orbifold cover is contractible and $n$ is even. Moreover, its orbifold fundamental group is an $n$-dimensional, rational Poincar\'e duality group.
\end{cor}

For a developable orbifold $M/G$ as above we can assume that $M$ is simply connected and that $G$ acts effectively; see Section \ref{sec:prel}. Then $M$ is the universal orbifold covering and $G$ is the orbifold fundamental group of $M/G$. Roughly speaking, a rational Poincar\'e duality group is a group whose rational (co)homology groups satisfy Poincar\'e duality; see Section \ref{sub:orb_coho} for the precise definition.

Dragomir derived several further necessary conditions for a developable compact Riemannian orbifold $M/G$ to have no closed geodesic. For instance, the Riemannian manifold $M$ itself cannot have closed geodesics and the group $G$ has to have a finite exponent \cite[Property~5.4]{Dragomir2}; cf. \cite{Adyan,Burnside}. 
To summarize, a negative answer to the following question would imply the existence of a closed geodesic on every compact Riemannian orbifold; see \cite[Theorems~4.3 and 4.7]{Dragomir2} and \cite{Dragomir}. Note that the orbifold fundamental group of a compact orbifold is always finitely presented \cite[Corollary I.8.11]{BH99} and that $G$ cannot contain an element of order $2$ by \cite[Theorem~4.7 (vi)]{Dragomir2} so that its exponent is necessarily odd. Further recall that an isometry is called elliptic if it has a fixed point.

\begin{qst}\label{qst:existence}
Is there a contractible, even-dimensional Riemannian manifold admitting a proper, cocompact action by elliptic isometries of a finitely presented group of finite odd exponent, such that the points with nontrivial isotropy are isolated?
\end{qst}

Such a Riemannian manifold has to be complete by Lemma \ref{lem:completeness} and the acting group must be infinite. It is not even known whether an infinite, finitely presented group of finite exponent exists. By Corollary \ref{cor:poincare_duality}, Question \ref{qst:existence} can only have a positive answer if there exists an infinite finitely presented rational Poincar\'e duality group of finite odd exponent. For further group theoretical considerations based on the non-existence of a closed geodesic on a developable Riemannian orbifold we refer to \cite[Chapter~5]{Dragomir2}.
\newline
\newline
\emph{Acknowledgements.} The authors thank Alexander Lytchak for stimulating discussions. The first named author moreover would like to thank Gudlaugur Thorbergsson in whose beautiful lectures he learned about the closed geodesic problem and the method to prove Theorem \ref{thm:B}. He was partially supported by the DFG funded project SFB/TRR 191. The support is gratefully acknowledged.

\section{Preliminaries}\label{sec:prel}

\subsection{Riemannian orbifolds}\label{pre:orbifolds}

For a definition of a Riemannian orbifold as a metric space we refer to \cite{La2}. The metric quotient $\Orb=M/G$ of a Riemannian manifold $M$ by a proper, isometric action of a discrete group $G$ is always a Riemannian orbifold; see e.g. \cite{Da11,La2}. A Riemannian orbifold $\Orb$ is called \emph{developable} or \emph{good} if it is isometric to such a quotient, in which case we refer to $M$ as a manifold cover of $\Orb$. If $\Orb$ is in addition connected, then, by orbifold covering space theory \cite{Th79}, it is isometric to such a quotient with a simply connected Riemannian manifold and an effective action. This simply connected manifold, which is uniquely determined by $\Orb$, is also referred to as universal orbifold covering of $\Orb$ and the effectively acting (deck transformation) group is the orbifold fundamental group of $\Orb$. The orbifold $\Orb$ is a manifold if and only if the deck transformation action on the universal orbifold cover is free.

\subsection{Orbifold geodesics}\label{sec:orbifold_geodesics}

 An (\emph{orbifold) geodesic} on a Riemannian orbifold is a path that locally lifts to a geodesic on a manifold chart; see \cite{La2}. In the developable case it can be described as a path that lifts to a geodesic on a manifold cover. A \emph{closed geodesic} is a nonconstant loop that is a geodesic on each subinterval. A closed geodesic on a developable Riemannian orbifold $M/G$ can be equivalently described as a nonconstant geodesic path $c: [0,1] \To M$ with the property that there exists some $g \in G$ such that $g(c(0))=c(1)$ and $dg(c'(0))=c'(1)$. If $M/G$ is compact, then $M$ is complete by Lemma \ref{lem:completeness} below and in this case such a path extends to a $g$-invariant geodesic on $M$. We see that for any subgroup $H$ of $G$ a closed geodesic on the Riemannian orbifold $M/H$ projects to a closed geodesic on $M/G$. In dimension two every compact Riemannian orbifold is either nondevelopable or finitely covered by a compact Riemannian manifold \cite{Sc83} and hence has a closed geodesic. In fact, there are  infinitely many closed geodesics in this case \cite{La1}.

\section{Closed geodesics on compact odd-dimensional Riemannian orbifolds}\label{sec:closed_geodesics}

A key ingredient for the induction step in the proof of Theorem \ref{thm:A} is the following statement.

\begin{lem}\label{lem:nor_coc} Suppose a discrete group $G$ acts properly, cocompactly and isometrically on a Riemannian manifold $M$. For $p\in M$ let $N$ be the fixed point set of the isotropy group $G_{p}$ of ${p}$, and let $H$ be the normalizer of $G_{p}$ in $G$. Then $H$ acts properly, cocompactly, and isometrically on $N$. 
\end{lem}
Note that the fixed point set $N$ in the statement of the lemma is a disjoint union of closed connected totally geodesic submanifolds all of which are open in $N$ but not necessarily of the same dimension, and that it is preserved by $H$. While $H$ does not in general preserve the connected components of $N$, it preserves the union of all components of $N$ of the same dimension which is a Riemannian manifold. Also note that the action of $H$ on $N$ does not need to be effective.

\begin{proof} The claim that the action of $H$ on $N$ is proper and isometric follows from the facts that it is a restriction of the proper and isometric action of $G$ on $M$ and that $H\times N$ is closed in $G \times M$. That the action is proper and isometric implies that the quotient topology on $N/H$ is induced by the quotient metric on $N/H$.

To show that the action of $H$ on $N$ is also cocompact, we first note that properness of the action implies that all isotropy groups $G_x$ are finite, and that for each $x$ there exists a neighborhood $U_x$ of $x$ such that $gU_x \cap U_x = \emptyset$ for all $g \in G \backslash G_x$; see \cite[Proposition~1.1]{Dragomir2}, for instance. In particular, all isotropy groups $G_y$ of points $y \in U_x$ are contained in $G_x$.

Suppose that the metric space $N/H$ is not compact. Then there exists a sequence of points $p_n \in N$ whose images in $N/H$ do not accumulate. By compactness of $M/G$, after passing to a subsequence there is a sequence $g_n \in G$ such that the sequence $x_n=g_np_n$ converges to a point $x \in M$. By passing to a subsequence we can assume that the isotropy groups $G_{x_n}$ do not depend on $n$ but are a fixed subgroup $\Gamma$ of $G_x$. Moreover, $G_{p}$ is contained in each $G_{p_n}$. Because of $G_{x_n}=g_n G_{p_n}g_n^{-1}$ this implies that for each $n$ the conjugated group $g_n G_{p} g_n^{-1}$ is contained in $\Gamma$. Again, by passing to a subsequence, we can assume that $g_n G_{p} g_n^{-1}$ does not depend on $n$. Therefore, each $g_n$ differs from $g_1$ by an element $h_n$ in the normalizer $H$ of $G_{p}$ in $G$, i.e. $g_n=g_1h_n$. It follows that $h_n p_n \in N$ converges to $y=g_1^{-1}x$. However, then the images of $p_n$ in $N/H$ converge to the image of $y$ in $N/H$ in contradiction to our assumption. This completes the proof of the lemma.
\end{proof}

More precisely, we need the following consequence of Lemma \ref{lem:nor_coc}.

\begin{cor}\label{cor:proper_coco_iso}
In the situation of Lemma \ref{lem:nor_coc} let $N_0$ be the connected component of $N$ containing the point $p$. Then $H_0=\{h\in H \mid h(N_0) = N_0\}$ acts properly, cocompactly, and isometrically on $N_0$.
\end{cor}
\begin{proof} The statement that the action of $H_0$ on $N_0$ is proper and isometric follows as in the proof of Lemma \ref{lem:nor_coc}.

To show cocompactness of the action we first note that the orbit map $\pi: N \rightarrow N/H$ is continuous and open. The connected component $N_0$ of $N$ is a closed subset of $N$. Since also the orbit $HN_0$ is a closed subset of $N$ as the complement of the union of the  (open) remaining connected components, we see that $\pi(N_0)$ is closed in $N/H$ and hence compact by Lemma \ref{lem:nor_coc}. The restriction of the orbit map $\pi$ to the connected component $N_0$, which is an open subset of $N$, is still continuous and open, and it factors through a continuous open map $\bar \pi: N_0/H_0 \rightarrow N/H$. We claim that $\bar \pi$ is injective. Indeed, for $x,y\in N_0$ with $\pi(x)=\pi(y)$ there is some $h\in H$ with $hx=y$ which implies that $h(N_0)=N_0$ by connectedness of $N_0$ so that in fact $h\in H_0$. Hence, we have a homeomorphism $\bar \pi: N_0/H_0 \rightarrow \pi(N_0)$ with compact image so that also $N_0/H_0$ is compact as claimed.
\end{proof}

We can now prove Theorem \ref{thm:A}.

\begin{proof}[Proof of Theorem \ref{thm:A}] Let $\Orb$ be an odd-dimensional compact Riemannian orbifold. To prove the claim, we can assume that $\Orb$ is connected. We can moreover assume that it is developable, i.e. of the form $M/G$ as above with a simply connected, hence orientable, Riemannian manifold $M$, as the existence of a closed geodesic on a nondevelopable orbifold is already known \cite{GH06}. After possibly passing to the orientation-preserving subgroup $G_+$ in $G$ of index at most $2$, we can moreover assume that $G$ preserves an orientation of $M$. Indeed, $M/G_+$ is an odd-dimensional compact Riemannian orbifold and the natural projection $M/G_+ \rightarrow M/G$ maps closed orbifold geodesics to closed orbifold geodesics. In particular, in the one-dimensional case, we can assume that $\Orb$ is a circle.

The rest of the proof works by induction on the dimension $n$ of $\Orb$. If $\Orb$ is a manifold, which it is for $n=1$, then the claim follows from the existence of closed geodesics on compact Riemannian manifolds \cite{LuFet}. Otherwise, we can choose a point $p \in M$ that projects to a maximal-dimensional singular stratum in $\Orb$, in other words, whose stabilizer group $G_p$ is nontrivial and for which the dimension of the fixed point set of the action of $G_p$ on $T_pM$ is maximal. Let $N_0$ be the connected component of the fixed point set of $G_p$ containing $p$. It is a closed totally geodesic submanifold of $M$ as discussed before. Since $G_p$ is nontrivial and since an isometry of a connected Riemannian manifold is the identity if its differential at some point is the identity, we see that $\dim N_0 < n$. By our maximality assumption the action of $G_p$ on the unit sphere in the normal space $T_pN_0^{\bot}\subset T_p M$ is free. Indeed, otherwise there is a nonzero $v \in (T_pN_0)^\bot$ fixed by some nontrivial $g\in G_p$ and for sufficiently small $\varepsilon>0$ the stabilizer $G_q$ of the point $q=\exp_p(\varepsilon v)$ is nontrivial as it contains $g$, and contained in $G_p$ as seen in the proof of Lemma~\ref{lem:nor_coc}. We can assume that $B_{2\varepsilon}(q)$ is a normal ball, in which case for each $\tilde q \in B_{2\varepsilon}(q) \cap N_0$ there is a unique minimizing geodesic from $q$ to $\tilde q$ whose endpoints are fixed by $G_q$ so that also its initial direction at $q$ is fixed by $G_q$. Since the minimizing geodesic between $q$ and $p$ meets $N_0$ transversally, the dimension of the span of all these directions is larger than the dimension of $N_0$ and so we obtain a contradiction to our maximality assumption. Since $G_p$ also preserves the orientation and since every orientation-preserving isometry of an even-dimensional sphere has a fixed point, we further see that the dimension of $N_0$ has to be odd, too. Let $H_0$ be the subgroup of the normalizer of $G_p$ in $G$ which preserves $N_0$. By Corollary~\ref{cor:proper_coco_iso} the quotient $N_0/H_0$ is an odd-dimensional, developable, compact Riemannian orbifold. Since $N_0$ is a totally geodesic submanifold of $M$, the natural map $N_0/H_0 \To M/G=\Orb$ maps closed geodesics to closed geodesics. As the dimension of $N_0$ is strictly smaller than the dimension of $M$, the claim hence follows by induction.
\end{proof}

\section{Closed geodesics on manifolds with cocompact isometric group actions}

\subsection{Closed geodesics on noncontractible Riemannian manifolds}

In the following we denote the path metric of a Riemannian manifold $M$ by $d$. For the proof of Theorem~\ref{thm:B} we first show the following fact.

\begin{lem}\label{lem:completeness} A connected Riemannian manifold $M$ with an isometric, cocompact group action is complete. Moreover, there exists a compact subset $K$ of $M$ whose translates under the action cover $M$ and the injectivity radius of $M$ is positive.
\end{lem}
\begin{proof} We denote the acting group by $G$ and consider the induced orbit pseudometric on $M/G$. Since its induced topology is coarser than the quotient topology, $M/G$ is also compact, and hence sequentially compact with respect to the pseudometric topology. For a Cauchy sequence $(p_n)$ in $M$ we can therefore, after passing to a subsequence, choose points $q_n$ in a $G$-orbit $X$ in $M$ with $d(q_n,p_n)<1/n$. In particular, $(q_n)$ is a Cauchy sequence as well. Let $B_r(q)$ be a normal ball in $M$ around some point $q \in X$ and let $g_n\in G$ be such that $g_nq_n=q$. We can choose $n_0\in \N$ such that $ g_{n_0}q_n \in \overline B_{r/2}(q)$ for all $n>n_0$. By compactness of $\overline B_{r/2}(q)$ the sequence $(g_{n_0}q_n)$ has a limit in $M$, and hence so does the sequence $(p_n)$.

To show the second claim we first observe that the diameter $D$ of the compact quotient pseudometric space $M/G$ is finite. Therefore, for any $p\in M$ the $G$-translates of the set $K=\overline B_{D+1}(p)$, which is compact by the Hopf--Rinow theorem, cover $M$. Positivity of the injectivity radius now follows from compactness of $K$ and invariance under the $G$-action.
\end{proof}

A classical approach to prove the existence of a closed geodesic on a compact Riemannian manifold $M$ is to choose a homotopically nontrivial map from a sphere into $M$, sweep out the sphere by loops and apply a min-max argument based on compactness and a curve shortening process in $M$, going back to ideas of Birkhoff \cite{Bi17,Bi27}, to find a sequence of loops that converges to a closed geodesic; see \cite{Kling,CM11}, for instance. In the noncompact case this strategy may fail as sequences of loops can leave every compact set. In the situation of Theorem \ref{thm:B} we can argue in a similar way and deal with the noncompactness by an application of Lemma \ref{lem:completeness}. In the following we explain more precisely how an elementary proof of the Lyusternik--Fet theorem based on  \cite[Appendix~A.1]{Kling}   generalizes to this setting, focusing on the necessary modifications and referring to \cite[Appendix~A.1]{Kling} for the standard technical details. The approach uses a discrete Birkhoff curve shortening-type process and requires no infinite-dimensional analysis.

\begin{proof}[Proof of Theorem \ref{thm:B}] 

We consider the space $\mathcal{P}$ of piecewise differentiable closed loops $c$ in $M$ defined on $[0,1]$ and equip it with the metric
$d_{\infty}(c_1,c_2)=\max_{t\in [0,1]} d\left( c_1(t),c_2(t) \right)$. The energy $E(c)=\frac{1}{2}\int_0^1\|\dot{c}(t)\|^2dt$ is defined on $\mathcal{P}$ and for $\kappa \geq 0$ we set $\mathcal{P}^{\kappa}=\{c\in \mathcal{P} \mid E(c) \leq \kappa\}$.

Since $M$ is connected and not contractible, for some $k\geq 0$ there exists a map $f:S^{k+1} \To M$ which is not null-homotopic by Whitehead's theorem \cite{Wh49}; see also \cite[p.~63]{Mi63}. By an approximation argument we can assume that the map $f$ is smooth.

The sphere $S^{k+1}$ is now swept out by a family of loops  $c_x:[0,1]\rightarrow S^{k+1}$ defined by $t\MTo   (x,\sqrt{1-\|x\|^2} \cos(2\pi t),\sqrt{1-\|x\|^2}\sin(2 \pi t))
$ for $x \in D^k$ and we set $f_x=f \circ c_x$. This induces a continuous map $F:(D^k,\partial D^k) \rightarrow (\mathcal{P}, \mathcal{P}^0)$ for which also the composition $E\circ F$ is continuous. By compactness we can assume that $\kappa>0$ is chosen in such a way that $F$ maps to $\mathcal{P}^\kappa$.

In the compact case there is a continuous energy-nonincreasing deformation $\mathcal{D}$ of $(\mathcal{P}^\kappa,\mathcal{P}^0)$ with the property that $E(\mathcal{D}(c))=E(c)$ for $c \in \mathcal P^\kappa$ if and only if $c$ is a closed geodesic \cite[Proposition~A.1.2]{Kling}. For its construction first some $m \in 2 \N$ with $4\kappa < \iota^2 m$ is chosen, where $\iota>0$ is the injectivity radius. Then every $c \in \mathcal{P}^{\kappa}$ satisfies 
$d\left(c(t),c(t+s)\right) <\iota$ for all $t$ and all $s\in [0,\frac{2}{m}]$ by the Cauchy--Schwarz inequality where the argument is understood modulo $1$. In particular, there exists a unique minimizing geodesic between $c(t)$ and $c(t+s)$.
The construction of the deformation now works by replacing arcs by minimizing geodesic segments in several steps; see \cite{Kling} for more details. We see that the construction of $\mathcal{D}$ extends to the setting of Theorem \ref{thm:B} since the injectivity radius is positive in this case by Lemma \ref{lem:completeness}. Moreover, since the construction is defined in terms of the geometry, it is equivariant with respect to isometries. Iterations of this deformation then give rise to maps $\mathcal{D}^nf:S^{k+1} \rightarrow M$ which are homotopic to $f$; see the proof of \cite[Theorem A.1.5]{Kling} for more details, also for the next sentence. Therefore, there is some $\varepsilon>0$, sufficiently small so that a loop of energy at most $\varepsilon$ is contained in a normal ball, such that for each $n\in N$ there is some $x\in \mathcal{D}^k$ with $E(\mathcal{D}^n f_{x})>\varepsilon$ as $f$ could be homotoped to a constant map otherwise. Choosing $x_n\in D^k$ such that $E(\mathcal{D}^nf_{x_n})=\max\{E(\mathcal{D}^nf_x)|x \in D^k\}:=e_n$ it then holds that $e_n \searrow e \ >0$.

Now let $K$ be a compact subset of $M$ as in the statement of Lemma \ref{lem:completeness} and let $K_{\kappa}$ be a closed $\sqrt{2\kappa}$-neighborhood of $K$ which is also compact. Again by the Cauchy--Schwarz inequality, for $n\geq 2$ we can choose $g_n \in G$ such that $g_n\mathcal{D}^{n-1}f_{x_n}=\mathcal{D}^{n-1}g_nf_{x_n}$ is contained in $K_{\kappa}$ and consider it as a geodesic $m$-gon. After passing to a subsequence we can assume by compactness that all $m$ vertices of $\mathcal{D}^{n-1}g_nf_{x_n}$ converge in $K_{\kappa}$. Then $\mathcal{D}^{n-1}g_nf_{x_n}$ converges uniformly to a geodesic $m$-gon $c_* \in \mathcal{P}^{\kappa}$; cf. \cite[Proposition A.1.1]{Kling}. Moreover, by continuity of the energy on the space of geodesic $m$-gons, the choice of the $x_n$, and continuity of $\mathcal{D}$ we have
\begin{equation*}
\begin{split}
E(\mathcal{D}c_*)&\leq E(c_*) = \lim_{n\To \infty} E(\mathcal{D}^{n-1}g_n f_{x_n})\leq \lim_{n\To \infty} e_{n-1}=e=E(\lim_{n\To \infty} \mathcal{D}^{n}g_n f_{x_n})= E(\mathcal{D}c_*)\,.
\end{split}
\end{equation*}
Hence, $E(\mathcal{D}c_*)=E(c_*)>0$, and so $c_*$ is a closed geodesic as discussed above.
\end{proof}

\subsection{Rational Poincar\'e duality groups}\label{sub:orb_coho} For the homological algebra background of the following discussion we refer to \cite{Brown}, for instance. A group $G$ is called an \emph{$n$-dimensional rational Poincar\'e duality group}, or $\mathrm{PD}^n_{\Q}$-group for short, if there exists a right $\Q G$-module $D$, isomorphic to $\Q$ as a $\Q$-module, and a (fundamental) homology class $\mu \in H_n(G;D)$ so that for any left $\Q G$-module $A$, cap product with $\mu$ defines an isomorphism $H^i(G;A)\cong H_{n-i}(G;D\otimes A)$ \cite{Bieri,DD89}. By the rational version \cite[Theorem~2.21, Definition~3.3]{DD89} of a result by Bieri--Eckmann \cite{BE73} (cf. \cite[Chapter VIII, Theorem 10.1]{Brown}) a group $G$ is a $\mathrm{PD}^n_{\Q}$-group if and only if it satisfies the following conditions 
\begin{enumerate}
\item $G$ is of type $FP_\Q$, i.e. $\Q$ admits a projective resolution over $\Q G$ of finite length by finitely generated $\Q G$-modules, and
\item $H^i(G;\Q G)\cong \begin{cases}
0, \text{ for } i\neq n\\
D, \text{ for } i= n
\end{cases}$
\end{enumerate}
with $D$ as above.

For the proof of Corollary \ref{cor:poincare_duality} we show the following lemma.

\begin{lem}\label{lem:pd_group} The orbifold fundamental group $G$ of a developable compact Riemannian $n$-orbifold $\Orb$ with contractible universal orbifold covering $M$ is a $\mathrm{PD}^n_{\Q}$-group.
\end{lem}
\begin{proof} We choose a $G$-equivariant triangulation $K$ of $M$ \cite{Il00}. By passing to a subdivision, for each $g\in G$ we can assume that a $g$-invariant simplex in $K$ is fixed pointwise by $g$. The $\Q G$-chain complex $C(K)$ has finite length and is finitely generated in each dimension by compactness of $\Orb$ and local finiteness of the chosen triangulation. Since the action of $G$ on $M$ has only finite isotropy groups, so does the action of $G$ on $K$. By an averaging argument as in the proof of Maschke's theorem the chain complex $C(K)$ is therefore projective; see \cite[Chapter~I, Exercise~8.5]{Brown}. Since $M$ is contractible, it gives rise to a projective resolution of $\Q$ over $\Q G$ satisfying condition $(i)$ above.

To verify condition (ii) above we first note that there is a natural isomorphism of cochain complexes of right $\Q G$-modules $\mathrm{Hom}_{\Q G}(C(K),\Q G) \cong \mathrm{Hom}_c(C(K),\Q)$, where $\mathrm{Hom}_c(C_i(K),\Q)$ consists of all $\Q$-homomorphisms $f:C_i(K)\rightarrow \Q$ such that, for every $a\in C_i(K)$, $f(ga)=0$ for all but finitely many $g\in G$. Indeed, the proof of the respective statement with $\Z$-coefficients \cite[Chapter~VIII, Lemma~7.4]{Brown} goes through with $\Q$-coefficients in the same way. The cohomology of $\mathrm{Hom}_{\Q G}(C(K),\Q G)$ is $H^*(G;\Q G)$, since $C(K)$ defines a projective resolution of $\Q$ over $\Q G$. Recall that $C(K)$ has a $\Q$-basis with one element for each simplex in $K$. These basis elements are permuted by $G$ with finite isotropy groups and fall into finitely many orbits. Therefore, $\mathrm{Hom}_c(C_i(K),\Q)$ consists of those cochains $f\in\mathrm{Hom}_c(C_i(K),\Q)$ with $f(\sigma)=0$ for all but finitely many simplices $\sigma$ in $K$. The cohomology of this complex is the cohomology of $K$ with compact support, denoted by $H_c^*(K;\Q)$; cf. \cite[Chapter~VIII, Proposition~7.5]{Brown}. Since $M$, as a contractible manifold, is orientable, Poincar\'e duality now implies that $H^*(G;\Q G)$
is concentrated in degree $n$ where it is isomorphic to $H_0(K;\Q)\cong\Q$ as a $\Q$-module. Hence, the claim follows with $D=H^n(G;\Q G)$.
\end{proof}

\subsection{Proof of Corollary \ref{cor:poincare_duality}} \label{sub:proof_cor}
Let $\Orb$ be a compact, connected Riemannian $n$-orbifold without a closed geodesic. By \cite{GH06}, $\Orb$ is developable, i.e. we can represent it as a quotient of a simply connected Riemannian manifold $M$, its universal orbifold covering, by a proper, cocompact, isometric action of a discrete group $G$; see Subsection \ref{pre:orbifolds}. Since the projection $M\rightarrow \Orb$ maps closed geodesics to closed geodesics (see Subsection \ref{sec:orbifold_geodesics}), $M$ has no closed geodesics either. Hence, $M$ is contractible by Theorem~\ref{thm:B}. By Theorem~\ref{thm:A} the dimension $n$ is even. The second part of the claim now follows from Lemma \ref{lem:pd_group}.


\begin{thebibliography}{[O'Ne66]}


\bibitem[Ad10]{Adyan} S. I. Adyan, \emph{The Burnside problem and related questions} (Russian); translated from Uspekhi Mat. Nauk 65 (2010), no. 5(395), 5--60 Russian Math. Surveys 65, no. 5, 805--855, (2010).

\bibitem[AS20]{Asselle} L. Asselle\ and\ F. Schmäschke, \emph{On geodesic flows with symmetries and closed magnetic geodesics on orbifolds}, Ergod. Theory Dyn. Syst., 40(6), 1480--1509, (2020).

\bibitem[BG92]{Benci} V. Benci\ and\ F. Giannoni, \emph{On the existence of closed geodesics on noncompact Riemannian manifolds}, Duke Math. J. {\bf 68}, no.~2, 195--215, (1992).

\bibitem[Bi72]{Bieri}  R. Bieri, \emph{Gruppen mit Poincar\'e-Dualität}, Comment. Math. Helv. 47, 373--396, (1972).

\bibitem[BE73]{BE73} R. Bieri and B. Eckmann, \emph{Groups with homological duality generalizing Poincar\'e duality}, Invent. Math. 20 (1973), 103--124. 

\bibitem[Bi17]{Bi17} 	G. D. Birkhoff, \emph{Dynamical systems with two degrees of freedom}, Trans. Amer. Math. Soc., 18 (1917) pp. 199--300.

\bibitem[Bi27]{Bi27} G. D. Birkhoff, \emph{Dynamical systems}, Amer. Math. Soc. Colloq. Publ., vol. IX, Amer. Math. Soc., Providence R.I., 1927.

\bibitem[B92]{Borzellino} J. E. Borzellino, {\it Riemannian geometry of orbifolds}, ProQuest LLC, Ann Arbor, MI, (1992).

\bibitem[BL96]{Borzellino2}  J. E. Borzellino\ and\ B. G. Lorica, \emph{The closed geodesic problem for compact Riemannian $2$-orbifolds}, Pacific J. Math. {\bf 175} (1996), no.~1, 39--46.

\bibitem[BH99]{BH99} M. R. Bridson and A. Haefliger, \emph{Metric spaces of non-positive curvature}, Grundlehren der Mathematischen Wissenschaften, vol. 319, Springer-Verlag, Berlin, 1999.

\bibitem[Bro94]{Brown}K. S. Brown, \emph{Cohomology of groups}, Graduate Texts in Mathematics, vol. 87, Springer-Verlag, 1994.

\bibitem[Bu02]{Burnside} W. Burnside, \emph{On an unsettled question in the theory of discontinuous groups}, Quart. J. Pure Appl. Math. 33 (1902), 230--238.

\bibitem[CM11]{CM11} T. H. Colding and W. P. Minicozzi II, A Course in Minimal Surfaces, Providence, R.I., American Mathematical Society, 2011.

\bibitem[Da11]{Da11} M. W. Davis, \emph{Lectures on orbifolds and reflection groups}, Transformation groups and moduli spaces of curves, Adv. Lect. Math. (ALM), vol. 16, Int. Press, Somerville, MA, 2011, pp. 63--93.

\bibitem[DD89]{DD89} W.~Dicks and M. J. Dunwoody, \emph{Groups acting on graphs}, Cambridge Stud. Adv. Math. 17, Cambridge University Press, Cambridge, 1989.

\bibitem[Dr11]{Dragomir2} G.~Dragomir, \emph{Closed geodesics on orbifolds}, PhD thesis, McMaster University, Hamilton, Ontario, (2011). 

\bibitem[Dr15]{Dragomir} G.~Dragomir, \emph{The Stratification of Singular Locus and Closed Geodesics on Orbifolds}, preprint, arXiv:1504.07157.

\bibitem[Gr74]{Gr74} K.~Grove, Isometry-invariant geodesics, Topology 13 (1974), 281–292.

\bibitem[Gr85]{Gr85} K.~Grove, The isometry-invariant geodesics problem: Closed and open. In: Geometry and Topology. Lecture Notes in Mathematics, vol 1167. Springer, Berlin, Heidelberg, (1985).

\bibitem[GT78]{GT78} K. Grove and M. Tanaka, On the number of invariant closed geodesics, Acta Math. 140 (1978), no. 1-2, 33–48.

\bibitem[GH06]{GH06} K. Guruprasad, A. Haefliger, \emph{Closed geodesics on orbifolds}, Topology, vol. 45 no. 3, 611--641, (2006).

\bibitem[Ha98]{Ha98}  J. Hadamard, “Les surfaces à courbures opposées et leurs lignes géodésiques,” Journal de Mathématiques Pures et Appliquées, Série 5, 4 (1898), 27--73.

\bibitem[Il00]{Il00} S. Illman, \emph{Existence and uniqueness of equivariant triangulations of smooth proper G-manifolds with some applications to equivariant Whitehead torsion},
J. Reine Angew. Math. 524 (2000), 129--183.

\bibitem[Kl78]{Kling} W. Klingenberg, {\it Lectures on closed geodesics}, Springer-Verlag, Berlin, 1978. 

\bibitem[Lan18]{La1} C. Lange, \emph{On the existence of closed geodesics on 2-orbifolds}, Pacific J. Math. {\bf 294} (2018), no.~2, 453--472.

\bibitem[Lan20]{La2} C. Lange, \emph{Orbifolds from a metric viewpoint}, Geom. Dedicata 209, 43--57 (2020).

\bibitem[LF51]{LuFet} L. A. Lyusternik\ and\ A. I. Fet, \emph{Variational problems on closed manifolds}, Doklady Akad. Nauk SSSR (N.S.) {\bf 81} (1951), 17--18.

\bibitem[Ma13]{Ma13} M. Mazzucchelli, On the multiplicity of isometry-invariant geodesics on product manifolds, Algebr. Geom. Topol. 14(1): 135-156 (2014).

\bibitem[Ma15]{Ma15} M. Mazzucchelli, Isometry-invariant geodesics and the fundamental group. Math. Ann. 362, 265–280 (2015).

\bibitem[Mi63]{Mi63} J. Milnor, \emph{Morse Theory}, Annals of Mathematics Studies 51, Princeton University Press, 1963.

\bibitem[P05]{Poincare} H. Poincar\'e, \emph{Sur les lignes g\'eod\'esiques des surfaces convexes},  Trans. Amer. Math. Soc. {\bf 6} (1905), 237--274.

\bibitem[Sc83]{Sc83} P. Scott, The geometries of 3-manifolds, Bull. London Math. Soc. 15 (1983), no. 5, 401–487.

\bibitem[Se01]{Secchi} S. Secchi, \emph{A note on closed geodesics for a class of noncompact Riemannian manifolds}, Adv. Nonlinear Stud. {\bf 1} (2001), no.~1, 132--142.

\bibitem[Th77]{Thorb1} G.~Thorbergsson,  \emph{Geschlossene Geodätische auf nichtkompakten Riemannschen Mannigfaltigkeiten}. Dissertation, Bonn, 1977.

\bibitem[Th78]{Thorb} G.~Thorbergsson, \emph{Closed geodesics on non-compact Riemannian manifolds}. Math. Z. 159 (1978), no. 3, 249--258. 

\bibitem[Th79]{Th79} W. P. Thurston, The geometry and topology of three-manifolds, lecture notes, Princeton University, 1979.

\bibitem[Wh49]{Wh49} J. H. C. Whitehead, \emph{Combinatorial Homotopy I}, Bulletin of the American Mathematical Society, 55(3), 213--245, (1949).
\end{thebibliography}
\end{document}